\documentclass[journal]{IEEEtran}
\ifCLASSINFOpdf
\else
\fi
\usepackage[pdftex]{graphicx}
\usepackage[outdir=./]{epstopdf}
\graphicspath{{./}}
\usepackage{caption}
\usepackage{subcaption}

\usepackage{amsmath} 
\usepackage{amssymb}  
\usepackage{mathtools}
\DeclarePairedDelimiter{\ceil}{\lceil}{\rceil}
\usepackage{cases}

\newtheorem{definition}{Definition}
\newtheorem{proposition}{Proposition}
\newtheorem{problem}{Problem}
\newtheorem{assumption}{Assumption}
\newtheorem{remark}{Remark}
\newtheorem{theorem}{Theorem}
\newtheorem{lemma}{Lemma}
\newtheorem{example}{Example}
\newtheorem{corollary}{Corollary}

\newcommand{\win}{\text{Win}}

\newcommand{\pre}{\text{Pre}}

\newcommand{\abs}[1]{\left\vert#1\right\vert}

\newcommand{\set}[1]{\left\{#1\right\}}
\newcommand{\Ball}[1]{\mathcal{B}_{#1}}

\newcommand{\sv}{\,\vert\;}
\newcommand{\st}{\text{ s.t. }}

\newcommand{\Real}{\mathbb{R}}
\newcommand{\Z}{\mathbb{Z}}

\newcommand{\U}{\mathcal{U}}
\newcommand{\D}{\mathcal{D}}
\newcommand{\X}{\mathcal{X}}

\renewcommand{\P}{\mathcal{P}}

\newcommand{\W}{\mathcal{W}}
\newcommand{\x}{\mathbf{x}}
\newcommand{\uu}{\mathbf u}
\newcommand{\dd}{\mathbf d}

\usepackage{enumerate}
\usepackage{enumitem}

\usepackage{algorithm}
\usepackage{algpseudocode}

\usepackage{tabularx}
\usepackage{multirow}
\usepackage{makecell}

\usepackage{xcolor}
\begin{document}
%
\title{\huge Robustly Complete Synthesis of Memoryless Controllers for Nonlinear Systems with Reach-and-Stay Specifications}
%
%
%

\author{Yinan~Li,~\IEEEmembership{Member,~IEEE,}
        Jun~Liu,~\IEEEmembership{Senior Member,~IEEE}
        \thanks{This work is supported in part by the NSERC DG, CRC, and ERA programs. Y. Li and J. Liu are with the Department of Applied Mathematics, University of Waterloo, Waterloo, Ontario, Canada. {\tt\small yinan.li, j.liu@uwaterloo.ca}}%
      }

\maketitle

\begin{abstract}
  This paper proposes {\color{black}a finitely terminating algorithm to solve reach-and-stay control problems for nonlinear systems. The algorithm is guaranteed to return a control strategy if the specification is robustly realizable.} Such a feature is desirable as the commonly used abstraction-based methods are sound but not complete for systems that are not incrementally stable. Fundamental to the proposed method is a fixed-point characterization of the \emph{winning set} of the system with respect to a given specification, i.e., the initial states that can be controlled to satisfy the specification. The use of an adaptive partitioning scheme not only guarantees the approximation precision of the winning set but also reduces computational time. The effectiveness and efficiency are illustrated by several benchmarking examples.
\end{abstract}

\section{Introduction}
\IEEEPARstart{R}{each}-and-stay control synthesis for nonlinear systems is concerned with finding control strategies that can steer the state of the system to a target set and maintain it in the target set afterwards. Such a problem exists in a variety of control applications such as voltage regulation of electrical power converters \cite{FribourgBook13}, attitude control and flight path following in flight control systems \cite{Faulwasser2009}, and regulation of room temperatures inside a building \cite{Nilsson2017}.

Seeking provably correct reach-and-stay control strategies at the presence of constraints and nonlinearity is challenging. Many nonlinear control methods, e.g., feedback linearization and Lyapunov-based control, are incapable to deal with constraints in system states or inputs. 
Integrating constraints and bounded perturbations in the stage of controller design was studied in \cite{Blanchini1992} for linear systems and extended to nonlinear systems \cite{Pin2014} where the reach-and-stay requirement is relaxed to reach a robustly controllable super set containing the target set when the target set is not controlled invariant. {\color{black}The set-theoretic approach used in \cite{Blanchini1992,Pin2014} is similar to earlier theoretical framework proposed in \cite{Bertsekas72} and \cite{Bertsekas71} for solving invariance and reachability control problems with linear or nonlinear dynamics and constraints. We note that the results in \cite{Bertsekas71,Blanchini1992,Pin2014} do not offer completeness guarantees, whereas the convergence in \cite{Bertsekas72} relies on exact computation of predecessor sets.} Model predictive control (MPC) was developed later as a standard framework for constrained control, but the feasibility of finding a controller is not guaranteed \cite{Grune2011}. 


More recently, abstraction-based methods \cite{tabuada2009verification}, which leverage model checking algorithms \cite{ClarkeBook1999} for control synthesis by constructing discrete abstractions of the original system dynamics, has gained popularity for solving various control problems. 
To be sound and complete in control synthesis, abstractions that are (approximately) equivalent to the original systems is usually needed, which is shown feasible for incrementally stable systems \cite{Pola2008,GirardPT10}. Without such a stability assumption, we can still construct over-approximations \cite{ZamaniPMT12,liu2013synthesis,liu2016finite,Reissig2016}, but it does not always guarantee a feasible control strategy, even if one exists, because spurious transitions are introduced and control synthesis is separated from abstraction. Using sufficiently small granularities, approximately complete control synthesis can be achieved without stability assumptions \cite{liu2017robust} but it is at the cost of intractable computation.

In this paper, we propose a \emph{sound and robustly complete} control synthesis method for discrete-time nonlinear systems to verify the existence of a reach-and-stay control strategy and construct the strategy if it exists. A control synthesis method is called robustly complete if it returns a control strategy whenever the given specification is realizable for the same system under bounded disturbances. This is an extension of our previous work for switched systems \cite{Li2018acc}. At the core is a fixed-point algorithm that determines the winning set over the continuous state space with respect to the given specification. Only assuming that the target set is compact, we prove that such an algorithm is sound and complete and a memoryless control strategy is sufficient for nonlinear systems while the one used in \cite{Blanchini1992} fails to be complete under the same condition.

One of our main contributions is that the proposed method is robustly complete, as opposed to most of the abstraction-based methods for nonlinear systems without stability assumptions (e.g., \cite{ZamaniPMT12,liu2013synthesis,Reissig2016,Nilsson2017,Hsu2018hscc}) that are not complete. Another benefit of our method is that it practically gains computational efficiency by partitioning only the region in the state space where necessary. This is achieved by an automatic subdivision framework based on interval computation \cite{jaulin2001applied}. In contrast with other applications of interval analysis such as \cite{Collins2008,Wan2009}, we deal with reach-and-stay control problems without the assumption that the target set is controlled invariant. Compared with the works with abstraction refinement mechanisms \cite{GirardGM16,Hsu2018hscc,Nilsson2017}, in which parameters need to be chosen empirically or synthesis does not always terminate in finite time, we device a scheme for adaptive tuning of discretization precision under a given threshold related to system robustness level.


\emph{\textbf{Notation}}: Let $\Z$, $\Z_{\geq 0}$, $\Real$, $\Real^n$ be the set of all integers, non-negative integers, reals and $n$-dimensional real vectors, respectively; $|\cdot|$ is the infinity norm in $\Real^n$ and $\textbf{1}^n$ indicates the $n$-dimensional vector with all elements equal to 1; given two sets $A,B\subseteq\Real^n$, $B\setminus A:=\set{x\in B\,\vert\, x\not\in A}$, $A\ominus B:=\{c\in\Real^n \,\vert\;c+b\in A, \forall b\in B\}$, and $\Ball{r}:=\set{y\in \Real^n \sv |y|\le r}$.

\section{Control Synthesis by Fixed-Point Iterations}

\subsection{The reach-and-stay control synthesis problem}
Consider nonlinear control system in the form of
\begin{align}
  \label{eq:df}
  x_{t+1}=f(x_t,u_t)+d_t,
\end{align}
where $f:\Real^n\times\Real^m\to\Real^n$, $ t\in \Z_{\geq 0}$, $x_t\in \X\subseteq\Real^n$ is the state, $u_t\in\U\subseteq\Real^m$ is the control input, $d_t\in \D\subseteq\Real^n$ is a bounded disturbance. 
The set $\X$ and $\U$ are assumed to be compact, and the bounded set $\D:=\set{d\in\Real^n\,|\,\abs{d}\leq \delta,\, \delta\geq 0}$. When $\delta=0$, system (\ref{eq:df}) reduces to the nominal one
\begin{align}
  \label{eq:df0}
  x_{t+1}=f(x_t,u_t).
\end{align}



A sequence of control inputs $\uu=\set{u_i}_{i=0}^\infty$, where $u_i\in\mathcal{U}$, is called a \emph{control signal}. Similarly, we denote by $\dd=\set{d_i}_{i=0}^\infty$ a \emph{disturbance}. A \emph{solution} of system (\ref{eq:df}) is an infinite sequence of states $\x=\set{x_i}_{i=0}^\infty$ generated by an initial condition $x_0\in\X$, a control signal $\uu$ and a disturbance $\dd$ such that (\ref{eq:df}) is satisfied for all $t\in\Z_{\geq 0}$.

\begin{definition}
  Let $\Omega$ be a subset of the state space $\X$ of a system (\ref{eq:df}). A \emph{reach-and-stay} property of a solution $\x=\set{x_i}_{i=0}^\infty$ of the system (\ref{eq:df}) with respect to $\Omega$, denoted by $\varphi(\Omega)$, requires that $x_k\in\Omega$ for all $k\geq j$, $k,j\in\Z_{\geq 0}$.
\end{definition}

The purpose of this paper is to design a control strategy, if there exists one, such that the resulting solution satisfies a given reach-and-stay objective $\varphi(\Omega)$. The set $\Omega$ will be omitted when the target area is clear from the context or we discuss a reach-and-stay objective in general. The form of control strategies is given in the following definition.

\begin{definition}
  A \emph{(memoryless) control strategy} of system (\ref{eq:df}) is a function $\kappa:\,\X\to 2^{\U}$.
A control signal $\uu=\set{u_k}_{k=0}^\infty$ is said to conform to a control strategy $\kappa$, if $u_k\in \kappa(x_k)$,$\forall k\ge 0$,
where $\set{x_k}_{k=0}^\infty$ is the resulting solution of system (\ref{eq:df}).
\end{definition}

If there exists an initial condition $x_0\in\X$ and a memoryless control strategy $\kappa$ such that, for any disturbance, the resulting solution of system (\ref{eq:df}) under a control signal that conforms to $\kappa$ satisfies the objective $\varphi$, we say $\varphi$ is \emph{realizable} for (\ref{eq:df}), and the control strategy $\kappa$ \emph{realizes} $\varphi$ for (\ref{eq:df}). The set of all realizable initial conditions is the \emph{winning set} of (\ref{eq:df}) with respect to $\varphi$, written as $\win^\delta(\varphi)$. For the nominal system (\ref{eq:df0}), the notation is simplified to $\win(\varphi)$. If $\win(\varphi)\neq \varnothing$ or $\win^\delta(\varphi)\neq \varnothing$, then $\varphi$ is realizable for system (\ref{eq:df0}) or (\ref{eq:df}).

\begin{problem}[Reach-and-Stay Control Synthesis]\label{pb:original}
  Consider a reach-and-stay objective $\varphi$ for system (\ref{eq:df}).
  \begin{enumerate}[label=(\roman*)]
  \item Determine if $\varphi$ is realizable, and
  \item synthesize a control strategy such that the closed-loop system satisfies $\varphi$ if possible.
  \end{enumerate}
\end{problem}

We assume neither that a target set $\Omega$ is controlled invariant (i.e., every state $x\in\Omega$ can be controlled inside $\Omega$ for all time) nor the existence of a control strategy, which is to be determined by solving Problem \ref{pb:original} (i).

\subsection{Fixed-point characterization of realizability}
The realizability of a reach-and-stay objective $\varphi(\Omega)$ with $\Omega\subseteq \X$ for system (\ref{eq:df}) can be determined through a fixed-point algorithm. The winning set of $\varphi(\Omega)$ is characterized by the fixed point, which is a subset of $\X$ returned by the algorithm. Considering that the system state space is always bounded in real applications, it is fair to assume the compactness of the state space $\X$ and the target set $\Omega$.

\begin{definition}\label{def:pre}
  Given a set $Y\subseteq\X$, the \emph{predecessor} of $Y$ with respect to system (\ref{eq:df}) is a set of states defined by
  \begin{align*}
    \pre^\delta(Y):=\{x\in \X \sv \exists u\in\U,\st f(x,u)+d\in Y,\forall d\in\D\}.
  \end{align*}
  The predecessor is simplified to $\pre(Y)$ for $\delta=0$. The set of \emph{valid control values} that lead to one-step transition to $Y$ for an $x\in\pre^\delta(Y)$ is 
  \begin{align*}
    U_Y(x):=\{u\in \U \sv f(x,u)+d\in Y,\, \forall d\in\D\}.
  \end{align*}
\end{definition}

For $X,Y\subseteq\X$, the predecessor $Y$ that resides in a set $X$ is the set $X\cap\pre^\delta(Y)$. To simplify the notation, we let 
\begin{align*}
  \pre^\delta(Y|X)=X\cap\pre^\delta(Y).
\end{align*}

The following properties are straightforward without further assumptions on system (\ref{eq:df}).
\begin{proposition}\label{prop:pre}
  Let $A, B\subseteq\X$ and $\delta\geq 0$. Then
  \begin{enumerate}[label=(\roman*)]
  \item $\pre^\delta(A)\subseteq \pre^\delta(B)$ if $A\subseteq B$,
  \item $\pre^\delta(X)=\pre(X\ominus\Ball{\delta})$ and $\pre^{\delta_2}(A)\subseteq \pre^{\delta_1}(A)$ if $0\leq\delta_1\leq\delta_2$.
  \end{enumerate}
\end{proposition}

If continuity is imposed to the dynamics, then $\pre^\delta(\cdot)$ has some additional property.
\begin{assumption}\label{asp:cont}
  The function $f:\Real^n\times\Real^m\to\Real^n$ in system (\ref{eq:df}) is continuous with respect to both arguments, and the state space $\X$ and the input space $\U$ are compact.
\end{assumption}

\begin{proposition}[\cite{rakovic2006}]\label{prop:closed}
  Under Assumption \ref{asp:cont}, if $A\subseteq\X$ is closed (compact), then $\pre^\delta(A)$ is closed (compact).
\end{proposition}

It is straightforward to see that Proposition \ref{prop:closed} still holds if $\U$ is a finite set. 
\begin{assumption}\label{asp:finite}
  The function $f:\Real^n\times\Real^m\to\Real^n$ in system (\ref{eq:df}) is continuous with respect to the first argument. The state space $\X$ is compact and input space $\U$ is finite.
\end{assumption}

We now present the following algorithm (\ref{eq:reachstay}) for reach-and-stay control synthesis, 
which consists of two nested fixed-point iterations.
\begin{align}
  \label{eq:reachstay}
  Y_\infty\Leftarrow
  \begin{cases}
    Y_0=\emptyset,X_0^\infty=\emptyset\\
    \begin{rcases}
      X_{i+1}^0=Y_i\cup\Omega\\
      X_{i+1}^{j+1}=\pre^\delta(X_{i+1}^j| X_{i+1}^j)
    \end{rcases}\Rightarrow X_{i+1}^\infty\\
    \kappa(x)\leftarrow U_{X_{i+1}^\infty}(x), \forall x\in \Omega\cap\left(X_{i+1}^\infty\setminus X_i^\infty\right)\\
    Y_{i+1}=\pre^\delta(X_{i+1}^\infty)\\
    \kappa(y)\leftarrow U_{X_{i+1}^\infty}(y), \forall y\in Y_{i+1}\setminus (Y_i\cup\Omega)
  \end{cases}
\end{align}




In the following proposition, we show that algorithm (\ref{eq:reachstay}) characterizes the winning set $\win(\varphi(\Omega))$. This implies that the realizability $\varphi(\Omega)$ for system (\ref{eq:df}) can be determined by checking the emptyness of $Y_\infty$. 

\begin{proposition}\label{prop:fp}
  Let $\Omega\subseteq\X$ be compact. Suppose that Assumption \ref{asp:cont} or \ref{asp:finite} holds. Let $Y_\infty=\bigcup_{i=0}^\infty Y_i$ be a fixed point of (\ref{eq:reachstay}), where $\set{Y_i}_{i=0}^\infty$ is a sequence of subsets of $\X$ generated from (\ref{eq:reachstay}). Then,
  \begin{enumerate}[label=(\roman*)]
  \item $\win^\delta(\varphi(\Omega))=Y_\infty$, and
  \item $\kappa$ is a memoryless control strategy that realizes $\varphi(\Omega)$.
  \end{enumerate}
\end{proposition}
\begin{IEEEproof}
  We only consider $\Omega\neq\emptyset$. Otherwise the results trivially hold. We first show $Y_\infty\subseteq\win^\delta(\varphi(\Omega))$ by induction. Trivially $Y_0=\emptyset\subseteq\win^\delta(\varphi(\Omega))$ and $X_1^\infty$ is compact. The induction step aims to show that, for all $i\geq 1$, $Y_{i+1}\subseteq\win^\delta(\varphi(\Omega))$ if $Y_i\subseteq\win^\delta(\varphi(\Omega))$. Assume that $X_i^\infty$ is compact. Then $Y_i=\pre^\delta(X_i^\infty)$ and thus $X_{i+1}^0=\Omega\cup Y_i$ is compact. The sequence $\{X_{i+1}^j\}_{j=0}^\infty$ is compact and decreasing by induction, using Proposition \ref{prop:pre} (i) and Proposition \ref{prop:closed} since $X_{i+1}^0=\Omega\cup Y_i$ is compact. It is also easy to show that $\set{Y_i}_{i=0}^\infty$ is increasing by induction. Furthermore, the compact limit set $X_{i+1}^\infty=\lim_{j\to\infty}X_{i+1}^j=\bigcap_{j=0}^\infty X_{i+1}^j$ (with respect to Painlev\'e-Kuratowski convergence \cite{RockafellarBook09}) is the maximal controlled invariant set inside $\Omega\cup Y_i$ \cite[Proposition 4]{Bertsekas72}. If $Y_i\subseteq\win^\delta(\varphi(\Omega))$, then $X_{i+1}^\infty\subseteq\win^\delta(\varphi(\Omega))$ because $X_{i+1}^\infty$ is a controlled invariant set inside $\Omega\cup Y_i$, which gives $Y_{i+1}=\pre^\delta(X_{i+1}^\infty)\subseteq\win^\delta(\varphi(\Omega))$ by Definition \ref{def:pre}. Hence, $\bigcup_{i=0}^\infty Y_i\subseteq\win^\delta(\varphi(\Omega))$.

  Applying the control inputs generated by $\kappa$, for all $i\geq 0$ and for all $x\in\Omega\cap\left(X_{i+1}^\infty\setminus X_i^\infty\right)$ and $x\in Y_{i+1}\setminus (Y_i\cup\Omega)$, the state $x$ will be controlled inside $\Omega\cup Y_i$ and $Y_i$ in one step, respectively. That means any state $x\in Y_{i+1}$ will be controlled into $Y_i$ until it enters $X_1^\infty\subseteq\Omega$, which is controlled invariant. Hence, we have also shown (ii) that $\kappa$ realizes $\varphi(\Omega)$.

  To see $\win^\delta(\varphi(\Omega))\subseteq Y_\infty$, we aim to show that $x\notin\win^\delta(\varphi(\Omega))$ for all $x\notin Y_\infty$. Let $x\notin Y_\infty$ be arbitrary. 
  Since $Y_\infty$ is a fixed point of (\ref{eq:reachstay}), i.e., $Y_\infty=\pre^\delta(V)$, where $V$ is the maximal controlled invariant set inside $\Omega\cup Y_\infty$.
  Then $x\notin \pre^\delta(V)$, which means that for all $\set{u_i}_{i=0}^\infty$ there exists $k$ and $\set{d_i}_{i=0}^k$ such that the resulting sequence of (\ref{eq:df}) satisfies $x_k\notin(\Omega\cup Y_\infty)$. Since $x_k\notin Y_\infty$, we can show in the same manner that for all $\set{u_i}_{i=k}^\infty$ there exists $k'\geq k$ and $\set{d_i}_{i=k}^{k'}$ such that the $k'$th state $x_{k'}$ of the resulting solution satisfies $x_{k'}\notin(\Omega\cup Y_\infty)$. In this way, for all $\set{u_i}_{i=0}^\infty$, we can find an infinite sequence $\set{d_i}_{i=0}^\infty$ for any $x\notin Y_\infty$ so that the resulting solution of (\ref{eq:df}) goes outside of $\Omega$ infinitely often. 
  Hence, $x\notin\win^\delta(\varphi(\Omega))$, which shows $\win^\delta(\varphi(\Omega))\subseteq Y_\infty$. The proof is now complete.
\end{IEEEproof}

{\color{black}Proposition \ref{prop:fp}} is a generalized result for the reach-and-stay problem. A algorithm for solving the reach-and-stay problem was first proposed in \cite{Blanchini1992} (as shown in (\ref{eq:oldreachstay})), which relies on the assumption that the target set is compact and convex.
\begin{align}
  \label{eq:oldreachstay}
  \begin{cases}
    \begin{rcases}
      X_0=\Omega\\
      X_{i+1}=\pre^\delta(X_i| X_i)\\
    \end{rcases}\Rightarrow X_\infty\\
    \kappa(x)\leftarrow U_{X_\infty}(x), \forall x\in X_\infty\\
    \begin{rcases}      
      Z_0=X_\infty\\
      Z_{i+1}=\pre^\delta(Z_i)\\
      \kappa(z)\leftarrow U_{Z_{i+1}}(z), \forall z\in Z_{i+1}\setminus Z_i
    \end{rcases}\Rightarrow Z_\infty
  \end{cases}
\end{align}

As opposed to (\ref{eq:reachstay}), algorithm (\ref{eq:oldreachstay}) is composed of two sequential fixed-point iterations, which fails to yield the real winning set. This can be illustrated in the following example.

\begin{example}
  Consider a target set $\Omega=[-0.3,0.3]\cup[0.8, 1.1]$ and the dynamics $x_{t+1}=-x_t(x_t^2-2.05x_t+0.05)+u_t+d_t$,
  where $x_t\in[-0.65,1.1]$, $u_t\in\set{0,10}$ and $d_t\in [-5,5]\times 10^{-4}$ for $t\in\Z_{\geq 0}$. We obtain the real winning set $\win^\delta(\varphi(\Omega))=Y_\infty=[-0.6311,1.1]$ by using algorithm (\ref{eq:reachstay}) while algorithm (\ref{eq:oldreachstay}) only gives a subset {\color{black}$Z_\infty=[-0.6311,-0.6082)\cup(-0.6021,0.9914)\cup(1.0135,1.1]$}. Figure \ref{fig:exp1} illustrates the difference between these two different algorithms.
  \begin{figure}
    \centering
    \includegraphics[scale=0.7]{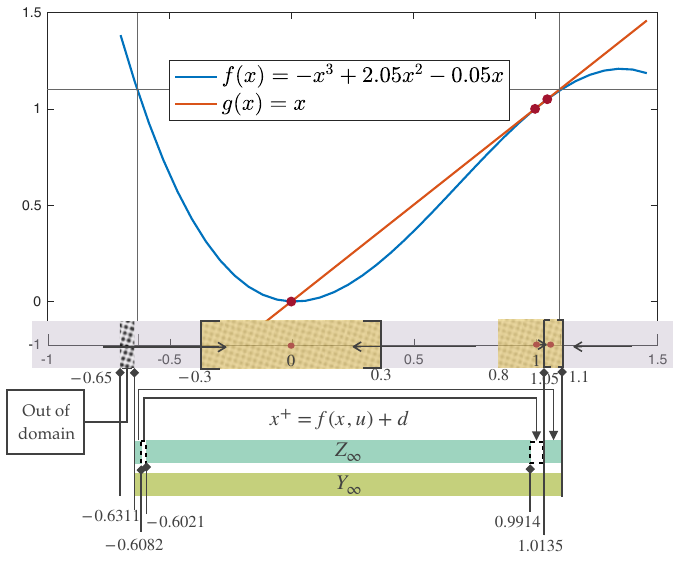}
    \caption{Let $u_t=0$ for all $t$. The fixed points $0$ and $1.05$ are stable while $1$ is unstable, which leads to $X_\infty=[-0.3,0.3]\cup[1.0370,1.1]$ if (\ref{eq:oldreachstay}) is used. Computing the reachable set of $X_\infty$ only gives a subset $Z_\infty$ of the real winning set.}
    \label{fig:exp1}
  \end{figure}
\end{example}


\begin{remark}
  In the literature, reach-avoid-stay objectives are also considered (e.g. \cite{Nilsson2017}), which additionally require the system state to avoid unsafe regions. By the definition of predecessor and algorithm (\ref{eq:reachstay}), the winning set, as well as every intermediate set, is bounded in the state space $\X$, which guarantees that any controlled trajectory using the synthesized control strategy $\kappa$ lives inside $\X$. Hence, (\ref{eq:reachstay}) can also be applied to solve reach-avoid-stay control problems by restricting the state space $\X$ to safe regions only.
\end{remark}



\subsection{Robustly complete control synthesis}
One problem with algorithm (\ref{eq:reachstay}), however, is that it might not terminate in a finite number of iterations.

\begin{example}
  Consider the system (in polar coordinates): $r_{t+1}= r_t^2$, $\theta_{t+1}=\mod(\theta_t+\theta_0/2\pi),\;\theta_0\in [0,2\pi)$.
  For this system, there is an unstable limit cycle given by $O=\set{(r,\theta)\in\Real\times[0,2\pi)\sv r=1}$.
Let the target set $\Omega$ be a subset of $O$ that contains the origin. The winning set $\win(\varphi)$ is the interior of $O$, which is open. Since the set returned by (\ref{eq:reachstay}) after a finite number of iterations is always closed, the algorithm cannot terminate in finite time.
\end{example}


Another difficulty is the computation of predecessors under nonlinear mappings. Only for some special cases, e.g. predecessors of polyhedral sets with respect to linear dynamics, which can be characterized by linear inequalities, the exact computation is possible. For general nonlinear dynamics, a practical way to overcome this difficulty is to use approximations. Inner-approximations are often used for control synthesis, since valid control values can not be found for all states in an outer-approximation. However, such a numerical compromise is at the expense of the loss of completeness, because the emptiness of the approximated winning set does not reflect the realizability of a given reach-and-stay property.

Alternatively, we study a relaxed problem based on the robust realizability of a specification introduced below. 

\begin{definition}
  A reach-and-stay objective $\varphi$ is said to be \emph{$\delta$-robustly realizable} for system (\ref{eq:df0}) if it is realizable for (\ref{eq:df}). If $\delta>0$, then $\varphi$ is called \emph{robustly realizable} for (\ref{eq:df0}).
\end{definition}

\begin{problem}[Robustly complete reach-and-stay control synthesis]\label{pb:relaxed}
  Consider a reach-and-stay property $\varphi$ for system (\ref{eq:df0}). Answer one of the two following questions:
  \begin{enumerate}[label=(\roman*)]
  \item Construct a control strategy if $\varphi$ is robustly realizable for system (\ref{eq:df0});
  \item Verify that $\varphi$ is not realizable for system (\ref{eq:df}) with $\delta>0$.
  \end{enumerate}
\end{problem}


\section{Reach-and-Stay Control Synthesis Is Robustly Complete} \label{sec:main}
This section presents our main results: there is a computer algorithm based on approximations of predecessors that will return reach-and-stay control strategies for nonlinear systems whenever the specification is robustly realizable. Such a conclusion is made because the inner-approximated winning set can be lower bounded by the winning set of the same system with additional perturbations, provided some condition to the precision of predecessor approximations is satisfied.

\subsection{An approximated synthesis for reach-and-stay control}
There are some computational redundancies in (\ref{eq:reachstay}): the sequence $\set{Y_i}_{i=0}^\infty$ is increasing and so is $\{X_i^j\}_{i=0}^\infty$ for all $j$. Hence, it is only necessary to compute the incremental parts between two adjacent sets in the sequences. Also, considering that predecessors cannot be precisely computed, we present the following approximated control synthesis algorithm (\ref{eq:interval}) based on an approximation $\widehat{\pre}$ of the predecessor operator $\pre$.
{\color{black}
\begin{align}
  \label{eq:interval}
  \widehat{Y}_\infty\Leftarrow
  \begin{cases}
    \widehat{Y}_0=\widehat{X}_0^\infty=\emptyset,V_0=\X\setminus\Omega\\
    \begin{rcases}
      W_{i}^0=\Omega\setminus \widehat{Y}_i\\
      \widehat{X}_{i+1}^j=\widehat{Y}_i\cup W_{i}^j\\
      W_{i}^{j+1}=\widehat{\pre}(\widehat{X}_{i+1}^{j}| W_{i}^j)
    \end{rcases}\Rightarrow
    \begin{aligned}
      W_{i}^\infty\\
      \widehat{X}_{i+1}^\infty
    \end{aligned}\\
    \kappa(x)\leftarrow U_{X_{i+1}^\infty}(x), \forall x\in W_{i}^\infty\\
    Z_{i}=\widehat{\pre}(\widehat{X}_{i+1}^\infty|V_i)\\
    \kappa(x)\leftarrow U_{X_{i+1}^\infty}(x), \forall x\in Z_{i}\\
    V_{i+1}=V_i\setminus Z_{i}\\
    \widehat{Y}_{i+1}=X_{i+1}^\infty\cup Z_{i}
  \end{cases}
\end{align}

\begin{theorem}\label{thm:main0}
  Let $Y_\infty$ and $Y_\infty^r$ ($r>0$) be the outputs of (\ref{eq:reachstay}) with operator $\pre$ and $\pre^r$, respectively. Suppose that $\widehat{\pre}(X)$ satisfies Propositions \ref{prop:pre}, \ref{prop:closed}, and
  $\pre^r(X)\subseteq\widehat{\pre}(X)\subseteq\pre(X)$ for any $X\subseteq\X$. Then 
  \begin{align*}
    Y_\infty^r\subseteq\widehat{Y}_\infty\subseteq Y_\infty.
  \end{align*}
\end{theorem}
}
\begin{IEEEproof}
  Let $\{\widetilde{X}_i^\infty\}$ ($\{\widetilde{X}_i^{r\infty}\}$) and $\{\widetilde{Y}_i\}$ ($\{\widetilde{Y}^r_i\}$) be the sequences of sets generated by algorithm (\ref{eq:interval}) with $\widehat{\pre}=\pre$ ($\widehat{\pre}=\pre^r$). We prove the theorem in the following structure: (i) Show that (\ref{eq:interval}) is equivalent to (\ref{eq:reachstay}) when set computation is accurate, i.e., $\widetilde{X}_i^\infty=X_i^\infty$ ($\widetilde{X}_i^{r\infty}=X_i^{r\infty}$) and $\widetilde{Y}_i=Y_i$ ($\widetilde{Y}^r_i=Y^r_i$) for all $i$. (ii) Show $\widetilde{Y}_\infty^r\subseteq\widehat{Y}_\infty\subseteq \widetilde{Y}_\infty$ under the given condition. 

  First of all, we show that $Y_i\subseteq\pre(Y_i)$ for all $i$. Since $X_i^\infty$ is a controlled invariant set, $X_i^\infty\subseteq\pre(X_i^\infty)$. By the definition of $Y_i$ in (\ref{eq:reachstay}) and monotonicity of $\pre$, $Y_i=\pre(X_i^\infty)\subseteq\pre(\pre(X_i^\infty))=\pre(Y_i)$. We now prove (i) by induction. The base case clearly holds since $\widetilde{Y}_0=Y_0=\widetilde{X}_0^\infty=X_0^\infty=\emptyset$. Suppose that $\widetilde{X}_i^\infty=X_i^\infty$ and $\widetilde{Y}_i=Y_i$ for some $i>0$. Then $\widetilde{X}_{i+1}^0=\widetilde{Y}_i\cup(\Omega\setminus\widetilde{Y}_i)=\widetilde{Y}_i\cup\Omega=X_{i+1}^0$, and
  \begin{align*}
    \widetilde{X}_{i+1}^{j+1}&=\widetilde{Y}_i\cup W_i^{j+1}=\widetilde{Y}_i\cup(\pre(\widetilde{X}_{i+1}^{j})\cap W_i^j)\\
                             &=(\widetilde{Y}_i\cup\pre(\widetilde{X}_{i+1}^{j}))\cap(\widetilde{Y}_i\cup W_i^j)\\
                             &=(\widetilde{Y}_i\cup\pre(\widetilde{X}_{i+1}^{j}))\cap\widetilde{X}_{i+1}^{j}.
  \end{align*}
  Also, $\pre(\widetilde{X}_{i+1}^{j})=\pre(\widetilde{Y}_i\cup W_i^j)\supseteq\pre(\widetilde{Y}_i)\supseteq \widetilde{Y}_i$, which implies that $\widetilde{X}_{i+1}^{j+1}=\pre(\widetilde{X}_{i+1}^{j})\cap\widetilde{X}_{i+1}^{j}$. This is the same as the iteration step in (\ref{eq:reachstay}), and thus $\widetilde{X}_{i+1}^\infty=X_{i+1}^\infty$. Now consider the sequence $\set{V_i}_{i=0}^\infty$. We have $V_0=\X\setminus\Omega$ and $V_{i+1}=V_i\setminus(\pre(\widetilde{X}_i^\infty)\cap V_i)=V_i\setminus\pre(\widetilde{X}_i^\infty)$. Unfolding $V_i$ until $V_0$ and using that $\pre(\widetilde{X}_i^\infty)\subseteq\pre(\widetilde{X}_{i+1}^\infty)$, we can derive $V_i=\X\setminus(\Omega\cup\pre(\widetilde{X}_i^\infty))=\X\setminus(\Omega\cup Y_i)=\X\setminus\widetilde{X}_{i+1}^0$. Then
  \begin{align}\label{eq:Y}
    \pre(\widetilde{X}_{i+1}^\infty)&=\pre(\widetilde{X}_{i+1}^\infty)\cap(\widetilde{X}_{i+1}^0\cup V_i)\nonumber\\
                                    &=\left[\pre(\widetilde{X}_{i+1}^\infty)\cap\widetilde{X}_{i+1}^0\right]\cup\left[\pre(\widetilde{X}_{i+1}^\infty)\cap V_i\right]\\
                                    &=\widetilde{X}_{i+1}^\infty\cup\pre(\widetilde{X}_{i+1}^\infty|V_i)=\widetilde{Y}_{i+1}.\nonumber
  \end{align}
  The equality $\pre(\widetilde{X}_{i+1}^\infty)\cap\widetilde{X}_{i+1}^0=\widetilde{X}_{i+1}^\infty$ can be derived by contradiction. If there exists $A\subseteq \widetilde{X}_{i+1}^0\setminus\widetilde{X}_{i+1}^\infty$ such that $A\subseteq\pre(\widetilde{X}_{i+1}^\infty)$ then $\widetilde{X}_{i+1}^\infty\cup A\subseteq\pre(\widetilde{X}_{i+1}^\infty\cup A)$, which indicates $A\cup\widetilde{X}_{i+1}^\infty$ is a larger controlled invariant set inside $\widetilde{X}_{i+1}^0$, but $\widetilde{X}_{i+1}^\infty$ is the maximal one. Therefore $Y_{i+1}=\widetilde{Y}_{i+1}$. The above argument also applies to prove $\widetilde{X}_i^{r\infty}=X_i^{r\infty}$ and $\widetilde{Y}^r_i=Y^r_i$.

  To prove (ii), we aim to show $X_i^{r\infty}\subseteq\widehat{X}_i^\infty\subseteq X_i^\infty$ and $Y_i^{r}\subseteq\widehat{Y}_i\subseteq Y_i$ for all $i$. Clearly $X^{r0}_1=\widehat{X}_1^0=X_1^0=\Omega$ and $\pre^r(X_1^{r0}|W_1^0)\subseteq\widehat{\pre}(\widehat{X}_1^0| W_1^0)\subseteq\pre(X_1^0|W_1^0)$ since $\pre^r(X)\subseteq\widehat{\pre}(X)\subseteq\pre(X)$ and Proposition \ref{prop:pre} (ii). This means $X^{r1}_1\subseteq\widehat{X}_1^1\subseteq X_1^1$. By induction, we can easily achieve $X^{rj}_1\subseteq\widehat{X}_1^j\subseteq X_1^j$ for any $j$. Thus $X_1^{r\infty}\subseteq\widehat{X}_1^\infty\subseteq X_1^\infty$. As shown in (\ref{eq:Y}), $\widehat{Y}_i=\widehat{\pre}(\widehat{X}_i^\infty)$. Then $Y^{r}_1=\pre^r(X_1^{r\infty})\subseteq\pre^r(\widehat{X}_1^{\infty})\subseteq\widehat{Y}_1\subseteq\pre(\widehat{X}_i^\infty)\subseteq\pre(X_i^\infty)=Y_1$. Therefore, (ii) can also be shown using induction.
\end{IEEEproof}

\subsection{Robustly complete control synthesis via interval arithmetic}
To implement the operator $\widehat{\pre}$ in (\ref{eq:interval}), we use interval arithmetic. This is because any compact set can be approximated by intervals with convergence guarantee under mild assumptions and interval operations are simple. An interval vector (box) in $\Real^n$ is denoted by $[x]$, where $[x]:=[x_1]\times\cdots\times[x_n]\subseteq \Real^n$ and $[x_i]=[\underline{x}_i,\overline{x}_i]\subseteq \Real$ for $i=1,\cdots,n$ with $\underline{x}_i$ as the infimum of $[x_i]$ and $\overline{x}_i$ the supremum. Let the set of all boxes in $\Real^n$ be $\mathbb{IR}^n$. The width of the interval $[x]$ is defined as $\text{wid}([x]):=\max_{1\leq i\leq n}\{\overline{x_i}-\underline{x_i}\}$. 

We have described in \cite[Algorithm 1]{Li2018rocs} a algorithm to obtain an inner approximation of $\pre(Y|X)$ with $\varepsilon$ precision ($\varepsilon>0$), denoted by $[\pre]^\varepsilon(Y|X)$ here, which is a union of a finite number of intervals. Central to this algorithm is the \emph{convergent inclusion function} $[f]:\,\mathbb{IR}^n\rightarrow\mathbb{IR}^m$ of $f:\,\Real^n\rightarrow\Real^m$ such that $f([x])\subseteq [f]([x])$ for all $[x]\in \mathbb{IR}^n$ and $\lim_{\text{wid}([x])\to 0}\text{wid}([f]([x]))=0$. 
The parameter $\varepsilon$ controls the minimum width of intervals for approximating $\pre(Y|X)$. We now discuss the relation between $\varepsilon$ and the error of set approximation in two scenarios.

\subsubsection{Finite control values}
In this case, system (\ref{eq:df}) can be treated as switched system, which has been discussed in \cite{Li2016,Li2018acc}. We summarize the result with its assumption as follows.

\begin{assumption} \label{asp:f}
  For system (\ref{eq:df}), there exists a constant $\rho>0$ such that
  \begin{align}\label{eq:lipschitz}
      |f(x,u)-f(y,u)|\leq \rho|x-y|,\quad \forall x,y\in \X.
  \end{align}
\end{assumption}

By Assumption \ref{asp:f}, we can always construct the centered-form convergent inclusion function $[f]([x],u)=f(\bar{x},u)+\rho([x]-\bar{x})\textbf{1}^n$ based on (\ref{eq:lipschitz}) for all $[x]\subseteq\X$.  

Under Assumption \ref{asp:f}, \cite[Lemma 1]{Li2016} can be used directly and presented as the following lemma.
\begin{lemma} \label{prop:cpre}
  Let $Y,X\subseteq \X$ be compact. If system (\ref{eq:df}) satisfies Assumption \ref{asp:f} in an neighborhood of $X$, then
  \begin{align*}
    \pre(Y\ominus\Ball{\rho\varepsilon}|X)\subseteq [\pre]^\varepsilon(Y|X)\subseteq \pre(Y|X).
  \end{align*}
\end{lemma}

\subsubsection{Infinite control values}
The compact set $\U\subseteq \Real^m$ might contain an infinite number of elements in $\U$. A straightforward way is to uniformly sample points in within the control set, e.g., an under-sampled set of controls
\begin{align}
  \label{eq:cset}
  [\U]_{\eta}:=\eta\Z^m\cap\U,
\end{align}
where $\Z^m$ denotes the $m$-dimensional integer grid, and $\eta\Z^m=\set{\eta z\sv z\in \Z^m, \eta>0}$. We additionally assume that for all $x\in\X$ and $u,v\in\U$,
\begin{equation}\label{eq:lipschitz3}
  |f(x,u)-f(x,v)|\leq \rho|u-v|.
\end{equation}


Similar to Lemma \ref{prop:cpre}, we prove the following approximation error by using under-sampled control values.

\begin{lemma}\label{prop:cpreu}
  Consider (\ref{eq:df}) with under-sampled control values (\ref{eq:cset}). Let $Y,X\subseteq\X$ be compact. If system (\ref{eq:df}) satisfies Assumption \ref{asp:f} and (\ref{eq:lipschitz3}) in a neighborhood of $X$, then
  \begin{align*}
    \pre(Y\ominus\Ball{\rho(\varepsilon+\eta)}|X)\subseteq [\pre]^\varepsilon(Y|X)\subseteq \pre(Y|X).
  \end{align*}
\end{lemma}
\begin{IEEEproof}
  We define a new predecessor operator $\pre_\eta(X):=\left\{x\in \X \sv \exists u\in[\U]_\eta,\st f_u(x)+d\in X,\forall d\in\D \right\}$.

  Let $Z=\pre(Y|X)$, $Z_\eta=\pre_\eta(Y|X)$ and $\widetilde{Y}=Y\ominus\Ball{\rho\frac{\eta}{2}}$. We first claim that $\pre(\widetilde{Y}|X)\subseteq Z_\eta\subseteq Z$. Trivially $Z_\eta\subseteq Z$ because $[\U]_\eta$ is a subset of $\U$. By Definition \ref{def:pre}, for all $z\in \pre(\widetilde{Y}|X)$, there exists a $u\in\U$ such that $f(z,u)+d\in\widetilde{Y}$ for all $d\in\D$. With (\ref{eq:lipschitz3}), for all $u\in\U$, there exists a $v\in[\U]_\eta$ such that $f(z,v)\in f(z,u)\oplus\Ball{\rho\frac{\eta}{2}}$. Then $f(z,v)+d\in f(z,u)\oplus\Ball{\rho\frac{\eta}{2}}+d=(f(z,u)+d)\oplus\Ball{\rho\frac{\eta}{2}}\in\widetilde{Y}\oplus\Ball{\rho\frac{\eta}{2}}=Y\ominus\Ball{\rho\frac{\eta}{2}}\oplus\Ball{\rho\frac{\eta}{2}}\in Y$ by \cite[Theorem 2.1 (ii)]{KolmanovskyG98}, which means that $z\in Z_\eta$. Hence the claim holds.

  By Lemma 1, $\pre_\eta(\widetilde{Y}\ominus\Ball{\rho\varepsilon}|X)\subseteq\underline{Z}\subseteq Z_\eta$. Applying the claim above, we have $\pre(\widetilde{Y}\ominus\Ball{\rho\varepsilon}\ominus\Ball{\rho\frac{\eta}{2}}|X)\subseteq \pre_\eta(\widetilde{Y}\ominus\Ball{\rho\varepsilon}|X)$. Therefore, $\pre(Y\ominus\Ball{\rho(\varepsilon+\eta)}|X)\subseteq \underline{Z}\subseteq Z_\eta\subseteq \pre(Y|X)$, which completes the proof.
\end{IEEEproof}

\begin{remark}\label{rmk:fxw}
  The above result can also be established for system $x_{t+1}=f(x_t,u_t,w_t)+d_t$, where $w_t\in \W\subseteq\Real^p$ and $\W$ is a bounded set of non-additive disturbances. As indicated in \cite[Lemma 1]{liu2017robust}, to achieve higher approximation precision in computing $[f]([x],u,\W)$, we can mince the set $\W$ into smaller sub-intervals and take the union of the images of all the sub-intervals under the inclusion functions. Suppose that $\W$ is uniformly partitioned with size $\mu$, i.e., $[\W]_\mu$. Replacing the inclusion function for (\ref{eq:df}) by $\bigcup_{[w]\in[\W]_\mu}[f]([x],u,[w])$, we can show, without much effort, that $\pre(Y\ominus\Ball{\rho(\varepsilon+\mu)}|X)\subseteq \underline{Z}\subseteq \pre(Y|X)$. 
\end{remark}

{\color{black}
\begin{theorem}\label{thm:main}
  Consider system (\ref{eq:df0}) under Assumption \ref{asp:finite}. Let $\Omega\subseteq\X$ be compact and Assumption \ref{asp:f} holds on $\X$. Suppose that $\varphi(\Omega)$ is $\delta$-robustly realizable for system (\ref{eq:df0}). Then algorithm (\ref{eq:interval}) with sets represented in intervals and $\widehat{\pre}=[\pre]^\varepsilon$ terminates in finite time and the following holds if $\rho\varepsilon\leq \delta$:
  \begin{align}
    \label{eq:cre}
    \win^\delta(\varphi(\Omega))\subseteq Y^\varepsilon \subseteq \win(\varphi(\Omega)).
  \end{align}
\end{theorem}
}
\begin{IEEEproof}
  Relation (\ref{eq:cre}) is an immediate result from Lemma \ref{prop:cpre} and Theorem \ref{thm:main0}, so we only show the finite termination here. Assume that $\Omega$ is an interval or a union of intervals. Under a given precision $\varepsilon>0$, $W_i$ can only be partitioned to finite number of intervals. Then, for the inner loop, there must exist a positive integer $N$ such that $W_i^N=\emptyset$ if $W_i^j\neq W_i^{j+1}$ for all $j\in\Z_{\geq 0}$ because $\{W_i^j\}_{j=0}^\infty$ is decreasing. Thus, the inner loop terminates within each outer loop. Likewise, the outer loop is also terminating since $\set{V_i}_{i=0}^\infty$ is decreasing and $\X$ only consists of a finite number of intervals.
\end{IEEEproof}

A similar result can also be established for systems under Assumption \ref{asp:cont}.

\begin{theorem}\label{thm:main2}
  Consider system (\ref{eq:df0}) under under Assumption \ref{asp:cont} with a set of under-sampled control values (\ref{eq:cset}). Let $\Omega\subseteq\X$ be compact and Assumption \ref{asp:f} and (\ref{eq:lipschitz3}) hold. Suppose that $\varphi(\Omega)$ is $\delta$-robustly realizable for system (\ref{eq:df0}). Then algorithm (\ref{eq:interval}) with sets represented in intervals and $\widehat{\pre}=[\pre]^\varepsilon$ terminates in finite time and the output $Y^\varepsilon$ satisfies (\ref{eq:cre}) if $\rho(\varepsilon+\eta)\leq \delta$.
\end{theorem}

\begin{remark}
  It is worth noting that the precision control parameter $\varepsilon$ in the inner and outer loops of algorithm (\ref{eq:interval}) can be set to different values, especially when the target area is volumetrically minuscule relative to the state space, e.g., in practical regulation problems. Furthermore, the precision control parameters are not necessarily fixed throughout the computation, but change with respect to the winning set obtained at each iteration.
\end{remark}

\begin{remark}
  The above theorems, however, cannot trivially lead to the convergence result, i.e., $\lim_{\varepsilon\to 0}Y^\varepsilon=\win(\varphi)$. This is because $\lim_{\delta\to 0}\win^\delta(\varphi)=\win(\varphi)$ does not always hold under Assumption \ref{asp:f}.
\end{remark}

By Theorems \ref{thm:main} (Theorem \ref{thm:main2}) and a controller defined in \cite[Proposition 3]{Li2018acc}, algorithm (\ref{eq:interval}) is guaranteed to generate a non-empty winning set along with a memoryless control strategy if the specification is robustly realizable. Even if algorithm (\ref{eq:interval}) returns an empty set, we can still make some conclusion on the robust realizability property of the specification, which are spelled out in the following corollary.

\begin{corollary}\label{co:main}
  Suppose that the assumptions for system (\ref{eq:df0}) in Theorem \ref{thm:main} (Theorem \ref{thm:main2}) hold. Then control synthesis with respect to a reach-and-stay specification $\varphi$ is robustly complete, i.e., there exists an algorithm that
  \begin{enumerate}[label=(\roman*)]
  \item generates a memoryless control strategy that realizes $\varphi$, if $\varphi$ is robustly realizable for system (\ref{eq:df0}), or
  \item verifies that $\varphi$ is not $\delta$-realizable for system (\ref{eq:df0}) for $\delta\geq\rho\varepsilon$ ($\delta\geq\rho(\varepsilon+\eta)$), if it returns no result. 
  \end{enumerate}
\end{corollary}

\begin{remark}\label{rmk:corollarymain}
  The conditions in Theorems \ref{thm:main} and \ref{thm:main2} serve as criteria for choosing the precision control parameter if the bound of disturbance $\delta$ and the Lipschitz constant $\rho$ over the state space can be trivially determined. Using such a criterion in actual computation is usually too conservative due to the evaluation of the Lipschitz constant over the entire state space. A practical benefit of Theorems \ref{thm:main} and \ref{thm:main2} is the guarantee that the winning set can be approximated more precisely by using a smaller precision parameter. 
Corollary \ref{co:main} implies that 
if we start computation with a large $\varepsilon$ and iteratively reducing it until the algorithm achieves a nonempty result, algorithm (\ref{eq:interval}) can also estimate the bound of the disturbances that can be tolerated without breaking the realizability of the given specification.
\end{remark}

\subsection{Example: automatic parallel parking}
We now demonstrate the effectiveness of the proposed algorithm on automatic parallel parking of the unicycle model \cite{AstromM08}: $\dot{x}=v \cos (\gamma+\theta) \cos (\gamma)^{-1}$, $\dot{y}=v \sin (\gamma+\theta) \cos (\gamma)^{-1}$, $\dot{\theta}=v \tan (\phi)$, where $(x,y)$ is the planar position of center of the unicycle, $\theta$ is its orientation, the control variable $v$ represents the velocity, $\phi$ is the steering angle command, and $\gamma=\arctan (a\tan (\phi)/b)$ with $a/b=1/2$.

In our simulation, the state space is $\X=[0,8]\times[0,4]\times[-72^\circ,72^\circ]$, sampling time is $\tau_s=0.3$s, and the set of control values is $\U=\set{\pm 0.9,\pm 0.6,\pm 0.3,0}$, which is sampled by uniform discretization of the space $[-1,1]\times[-1,1]$ with grid width $\eta=0.3$. An exact discrete-time model of the unicycle can be obtained \cite{Li2018arxiv} and readily verified Lipschitz continuous over $\X$ for all control values in $\U$.

Suppose that the length and width of the unicycle be $L=2$ and $H=1$, respectively. We consider two problem settings: parking with a wide marginal space $\Delta=L=2$ and a narrow marginal space $\Delta=0.5$. The marginal space is the distance between the front and rear vehicles in addition to $L$. For both cases, the rear vehicle center is at $(1,0.5)$, and thus the front vehicle center is at $(1+3L/2+\Delta, 0.5)$. The target area is $\Omega=[1+L, 1+L+\Delta]\times[0.5,0.6]\times[-3^\circ,3^\circ]$.

\begin{figure}[htbp]
  \centering
  \begin{subfigure}[htbp]{0.5\linewidth}
    \includegraphics[width=1\linewidth]{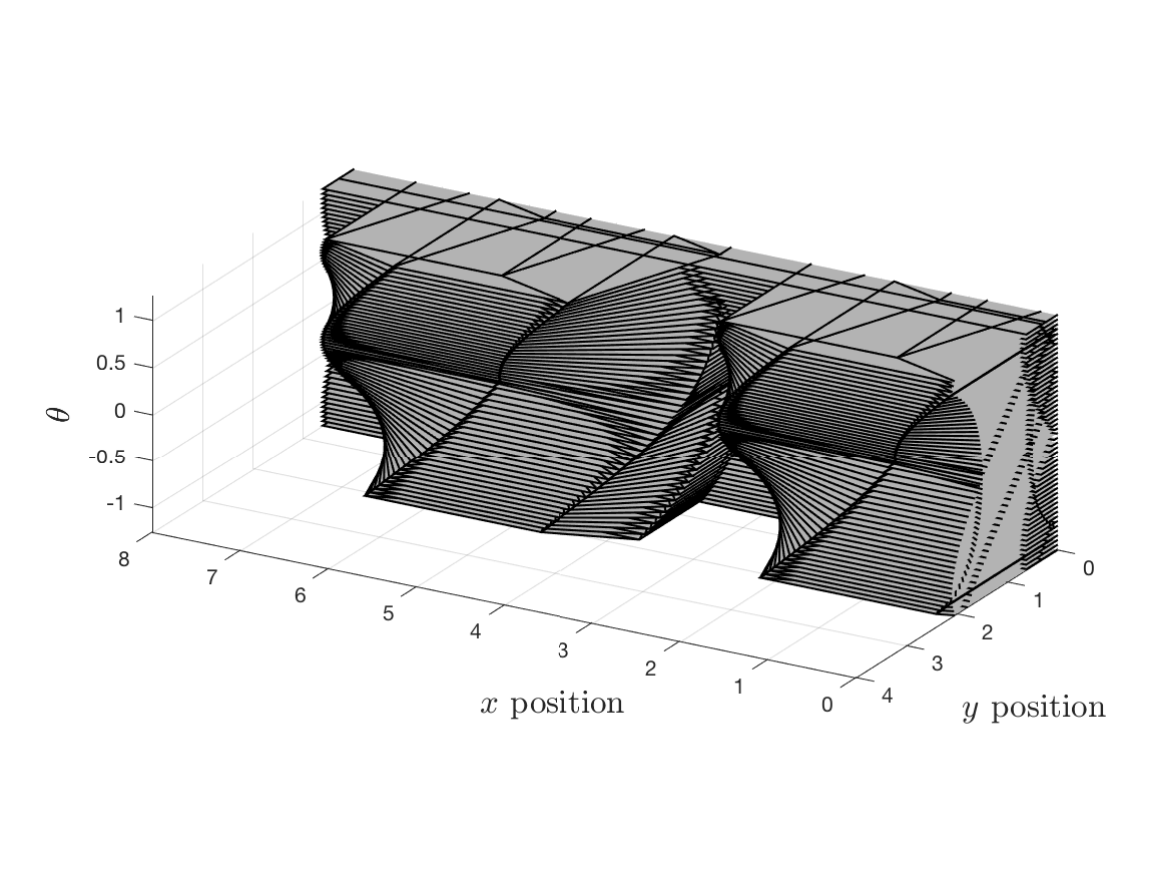}
    \caption{The $x-y-\theta$ view.}
  \end{subfigure}%
  \begin{subfigure}[htbp]{0.5\linewidth}
    \includegraphics[width=1\linewidth]{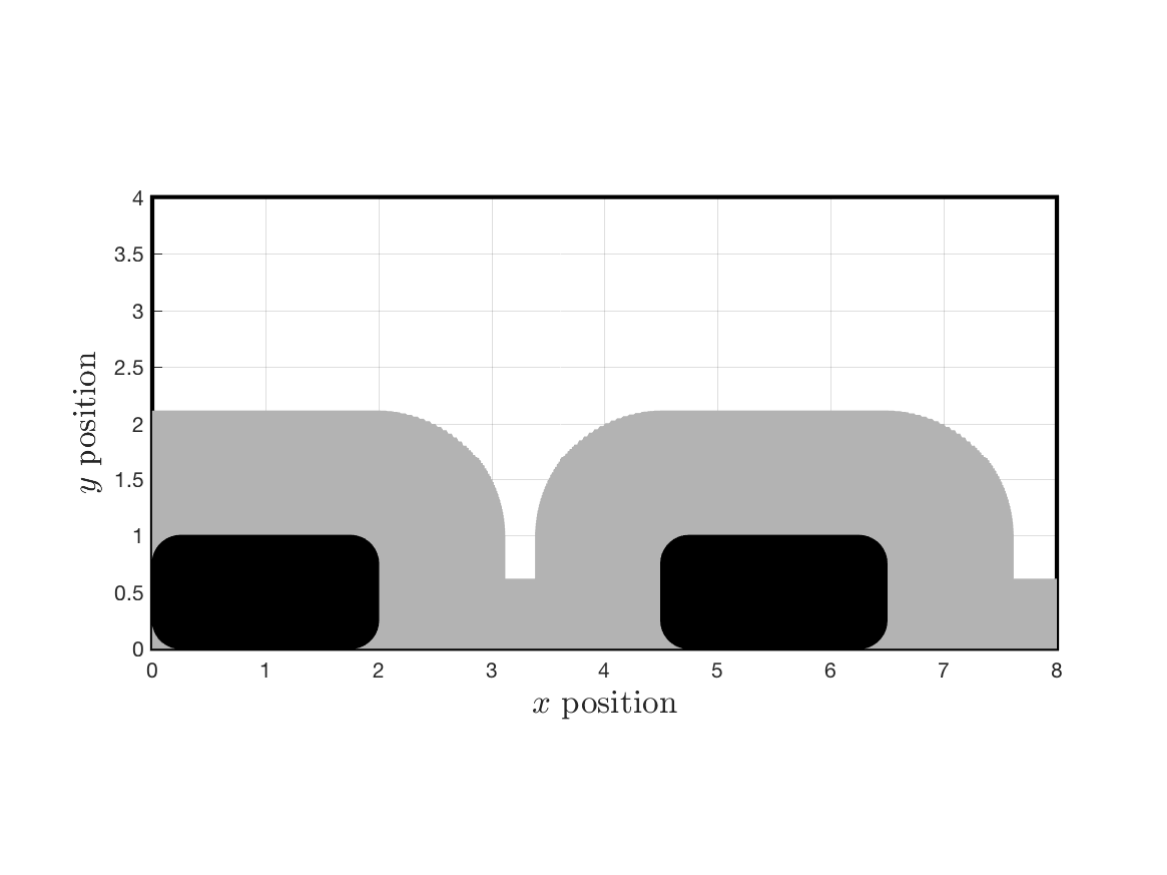}
    \caption{The $x-y$ view.}
  \end{subfigure}
  \caption{Collision area when $\Delta=0.5$. In (b), the gray area is the $x-y$ plane projection of the 3D collision area, and the two black rectangles represent the bodies of rear and front vehicle.}
  \label{fig:collisionarea}
\end{figure}

The collision area (the center position and orientation of the unicycle that causes collision with the parked vehicles and the curb) needs to be determined before control synthesis. We assume that vehicles and the curb are rectangles. Then the collision area can be interpreted by inequalities of the form $g(x)\leq 0$, which is derived by checking if two polyhedra intersect. It is clear that the center of the unicycle has different admissible regions with different orientations. Hence, the collision area is not simply a hyper-rectangle in $\Real^3$, as shown in Fig. \ref{fig:collisionarea} (a). The free configuration space (the admissible position of the unicycle center in $\Real^3$) determined by such a constraint can be handled by algorithm (\ref{eq:interval}).

We perform control synthesis for both cases using ROCS \cite{Li2018rocs}. 
By Corollary \ref{co:main} (i), if parallel parking is robustly realizable with the given marginal space, we can always synthesize a control strategy using a sufficiently small precision without calculating the Lipschitz constant. To see if the specifications in these two parking scenarios are realizable, we use different precision control parameters. The corresponding control synthesis results regarding the number of partitions (\#$\P_{1,2}$) and the run time ($t_{1,2}$) are summarized in TABLE \ref{tbl:parking}.
\begin{table}[htbp]
  \centering
  \caption{Control synthesis with different precisions.}
  \label{tbl:parking}
  \begin{tabular}{c||c|c||c|c} 
    \hline
    $\varepsilon$ & \#$\P_1$ & $t_1$ (s) & \#$\P_2$ & $t_2$ (s)\\
    \hline
    0.07 & 176786 & 102.93 & -- & -- \\
    \hline
    0.06 & 176666 & 103.19 & 1797027 & 295.68 \\
    \hline
    0.02 & 203166 & 127.44 & 1832589 & 327.50 \\
    \hline
    0.01 & 274694 & 176.20 & 1920929 & 427.48 \\
    \hline
  \end{tabular}
\end{table}

\begin{figure}[htbp]
  \centering
  \begin{subfigure}[htbp]{0.5\linewidth}
    \includegraphics[width=0.95\linewidth]{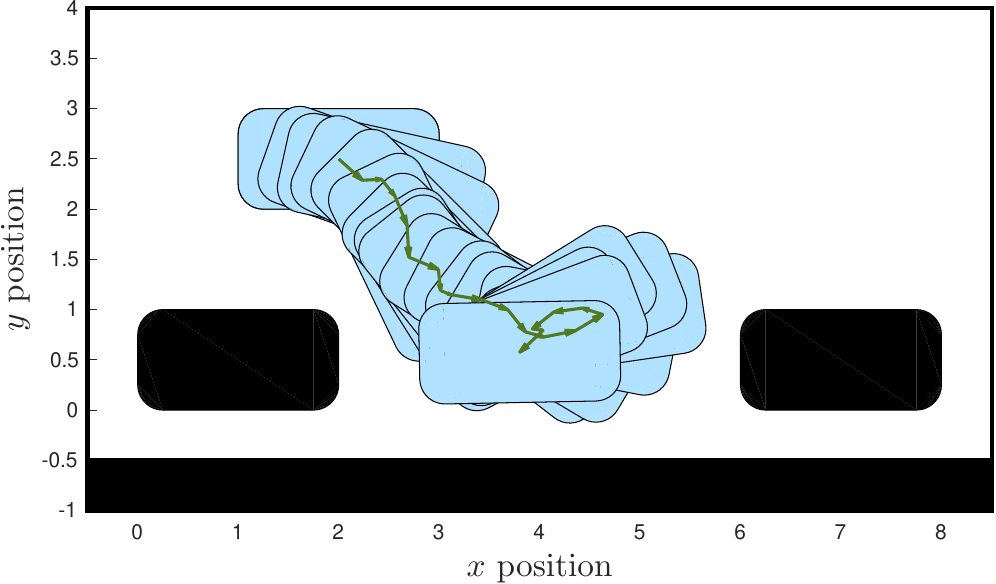}
    \caption{$\Delta=2$, $(x_0,y_0)=(2, 2.5)$.}
  \end{subfigure}%
  \begin{subfigure}[htbp]{0.5\linewidth}
    \includegraphics[width=0.95\linewidth]{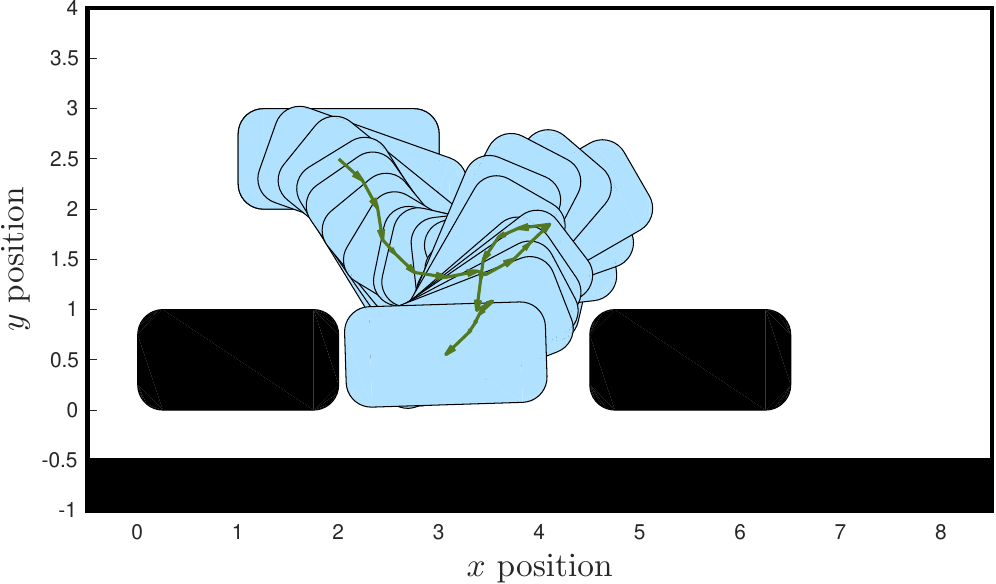}
    \caption{$\Delta=0.5$, $(x_0,y_0)=(2, 2.5)$.}
  \end{subfigure}
  \begin{subfigure}[htbp]{0.5\linewidth}
    \includegraphics[width=0.95\linewidth]{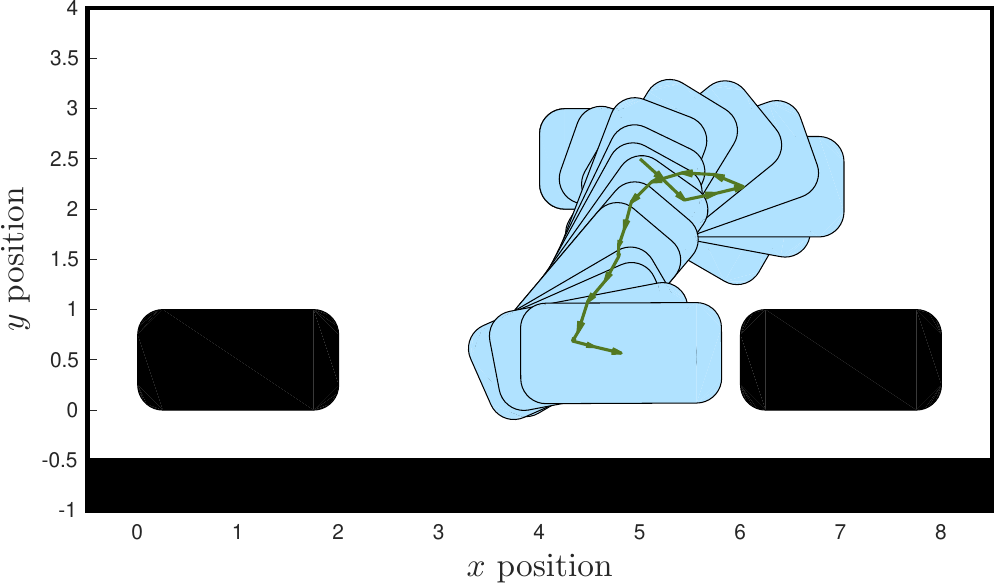}
    \caption{$\Delta=2$, $(x_0,y_0)=(5, 2.5)$.}
  \end{subfigure}%
  \begin{subfigure}[htbp]{0.5\linewidth}
    \includegraphics[width=0.95\linewidth]{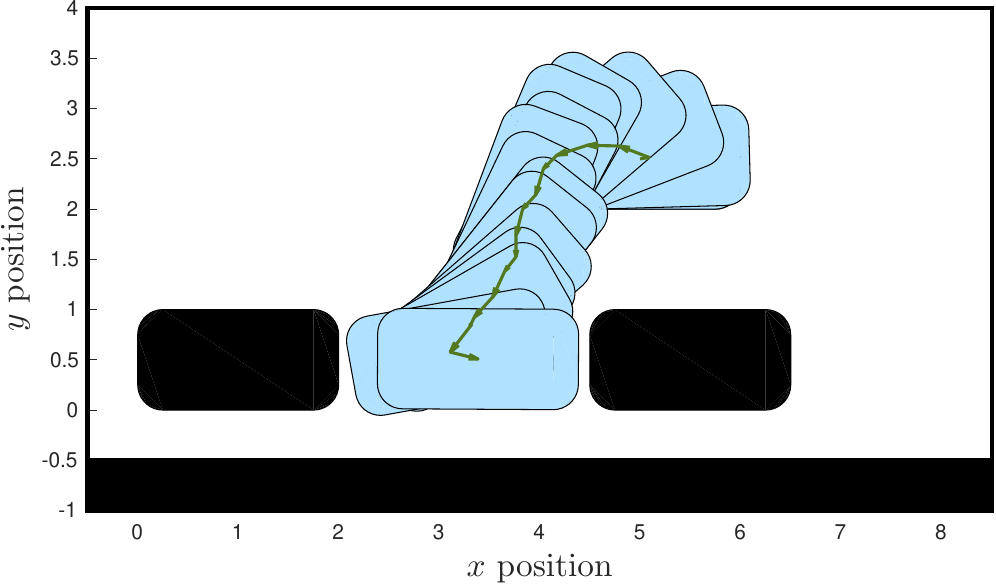}
    \caption{$\Delta=0.5$, $(x_0,y_0)=(5, 2.5)$.}
  \end{subfigure}
  \caption{Controlled parking trajectories from an initial condition $(x_0,y_0)$ with wide and narrow marginal parking spaces.}
  \label{fig:pparking}
\end{figure}

For both scenarios, the unicycle can be successfully parked into the target spot from any point of the free configuration space. The controlled parking trajectories with the resulting memoryless control strategies are presented in Fig. \ref{fig:pparking}, which all meet the parallel parking specification.

When the marginal parking space $\Delta$ is 0.5, we need a control synthesis precision no greater than 0.06 so that a memoryless control strategy can be generated. Additionally for this specific example, using a smaller $\varepsilon$ only increases the winning set by adding intervals close to the boundary of the free configuration space.

{\color{black}
\section{Evaluation of Time Complexity}
To show how well the proposed method performs in terms of computational time, we compare the time complexities of abstraction-based methods and the proposed method. Although theoretical analysis shows the equivalency of both methods in the worst case, the proposed method outperforms abstraction-based methods in solving many of the control synthesis problems practically.

\subsection{Complexity analysis}
Let $\varepsilon$ and $\eta$ ($\varepsilon,\eta>0$) be the grid size of the state space $\X$ and input space $\U$, respectively. Assume that the cost in terms of run time for each computation of the predecessor is some constant $c>0$, and $c_1,c_2>0$ are some constants related to the width of the state and input space.

For abstraction-based control synthesis based on a uniform partition of the state space, the number of discrete states and inputs are $N_S=\ceil{\frac{c_1}{\varepsilon}}^n$ and $N_U=\ceil{\frac{c_2}{\eta}}^m$, respectively. Then the time complexity for computing abstractions is $\mathcal{O}(cN_SN_U)$. Under Assumption \ref{asp:f}, the number of transitions $N_T$ is $(\ceil{\rho}+1)^nN_SN_U$. Using the classical co-B\"uchi algorithm, the time for solving the discrete control synthesis problem is $\mathcal{O}(N_SN_T)$, which yields the overall time complexity of abstraction-based control synthesis:
\begin{align}\label{eq:complexity-abst}
  \mathcal{O}(cN_SN_U+(\ceil{\rho}+1)^nN_S^2N_U).
\end{align}

We now analyze the time complexity of algorithm (\ref{eq:interval}) via interval computation implemented based on a binary tree data structure. According to the bisection scheme for predecessor approximation, the greatest depth of the binary tree is {\color{black}$h_{\max}=\ceil{n\log_2(\frac{c_1}{\varepsilon})}\approx\log_2 N_S$}. Set membership tests are performed by searching the binary tree. Hence, in the worst case where the tree is of depth $h_{\max}$, computation of each predecessor, including membership test, takes approximately $(h_{\max}+c)N_U$ operational time.
Let $N_G$ be the number of the set of intervals that represents the target set $\Omega$. Then the number of intervals outside of $\Omega$ is $N_S-N_G$. In the worst case for algorithm (\ref{eq:interval}), the set elements in the sequences $\set{Y_i}_{i=0}^\infty$ and $\{X_i^j\}_{j=0}^\infty$ differ by one interval. Then the number of iterations $N_I$ is
\begin{align*}
  \sum_{i=0}^{N_G}(i^2+i)+\sum_{i=0}^{N_S-N_G}i&=\frac{N_G^3+3N_G^2+8N_G}{6}\\
                                              &+\frac{(N_S-N_G)^2+(N_S-N_G)}{2}.
\end{align*}
If $N_G<<N_S$, then $N_I\approx(N_S^2+N_S)/2$. Hence, the time complexity of the algorithm (\ref{eq:interval}) is of
{\color{black}
\begin{align}
  \label{eq:complexity-interval}
  \mathcal{O}(\frac{c}{2}N_UN_S^2+\frac{1}{2}N_UN_S^2\log_2 N_S).
\end{align}

By comparing (\ref{eq:complexity-abst}) with (\ref{eq:complexity-interval}), the time complexity of algorithm (\ref{eq:interval}) is of $\mathcal{O}(N_UN_S^2\log N_S)$ while the abstraction-based methods is quadratic in $N_S$. 
}The overhead of algorithm (\ref{eq:interval}) primarily comes from the set inclusion tests by searching the binary tree, i.e., the part induced by $h_{\max}$. When only a high precision is necessary to yield a control strategy, the overhead run time is relatively large, which makes (\ref{eq:interval}) less efficient than abstraction-based methods. 

{\color{black}The worst case, however, rarely exists in practical control problems. On the other hand, the use of a non-uniform partitioning scheme avoids partitioning the region in the state space without helping in control synthesis. This usually leads to fewer discrete states for a given precision.} In this sense, the proposed method is less sensitive to the state and input discretization precisions than abstraction-based control synthesis methods. 
Such results from complexity analysis will be shown by the comparison tests in the following section.

From the relationship between system dimension and the time complexity as discussed above, the main limitation of the proposed method, which also exists in abstraction-based methods, is that it still suffers the \emph{curse of dimensionality}.

\begin{table*}[t]
  \centering
  \caption{Performance comparison tests: TO = \emph{time out} ($>86400$ s) and ``--'' = \emph{control synthesis fails}.}
  \label{tab:compare}
  \begin{tabular}{cc||c|c|c||c|c|c||c|c|c|c|c}
    \hline
    \multicolumn{2}{c||}{Examples} & \multicolumn{3}{c||}{Parameters} & \multicolumn{3}{c||}{ROCS} & \multicolumn{5}{c}{SCOTS}\\
    \hline
                                   & & \multirow{2}{*}{$n$} & \multirow{2}{*}{$N_U$} & \multirow{2}{*}{$\varepsilon$} & \multirow{2}{*}{$N_S$} & \multirow{2}{*}{\#Iter} & \multirow{2}{*}{time(s)} & \multirow{2}{*}{$N_S$} & \multirow{2}{*}{$N_T$} & \multirow{2}{*}{\#Iter} & \multicolumn{2}{c}{time(s)}\\
    \cline{12-13}
                                   & & & & & & & & & & &Abst & Syn\\
    \hline
    \multicolumn{2}{c||}{\multirow{2}{*}{DC-DC converter}} & \multirow{2}{*}{2} & \multirow{2}{*}{2} & 0.005 & 22433 & 76(529) &0.53 &40401 &291068 &84(671) &0.69 &15.90 \\
    \cline{5-13}
                                   & & & & 0.001 & 162261& 76(272)& 3.48& 1002001& 7243320 & 77(431)& 29.83 & 481.97\\
    \hline
    \multicolumn{2}{c||}{\multirow{2}{*}{Motion planning}}
                                   & \multirow{2}{*}{3} & \multirow{2}{*}{49} & 0.2 & 280291& 381(1)& 151.01& 91035& $3.73\times 10^7$&-- & 82.80&--\\
    \cline{5-13}
                                   & & & & 0.1 &1850830 &297(1) &1062.97 & 724271&$2.95\times 10^{8}$ & 313(2266) &2004.66 &17568.2 \\
    \hline
    \multirow{2}{*}{\makecell{Parallel\\ parking}} & $\Delta=0.5$& \multirow{2}{*}{3} & \multirow{2}{*}{49} & 0.02& 1832589& 133(1)& 327.50& 10075125& TO & TO & TO & TO\\
    \cline{5-13}
                                   & $\Delta=2$&  &  & 0.07& 167155& 123(8)& 94.32& 83025& $3.277\times 10^7$&-- & 73.14&--\\
    \hline
  \end{tabular}
\end{table*}


\subsection{Experimental tests on performance}
We now compare the performance of our proposed method with abstraction-based methods (implemented in SCOTS \cite{scots}) on solving different benchmarking examples. The results on a 3.6 GHz processor (Intel Core i3) are shown in TABLE \ref{tab:compare}. The column of \#Iter indicates the number of outer loops and the total number of inner loops (in the bracket) running (\ref{eq:reachstay}). The time for abstraction-based control methods is split into the part for abstraction (indicated as \emph{Abst}) and the one for synthesis (as \emph{Syn}).

The state space in the \emph{DC-DC converter} example (see \cite{Li2016} for the detailed model) is $\X=[0.649,1.65]\times[0.9898,1.19]$ and the target region $\Omega$ of the reach-and-stay specification $\varphi(\Omega)$ is $[1.1,1.6]\times[1.08,1.18]$. While the full setting of the \emph{Motion planning} example can be found in multiple works (e.g. \cite{ZamaniPMT12,Reissig2016}), we consider in our experiments the reach-and-stay control objective instead of just reachability. 

Our proposed method outperforms abstraction-based methods in those examples. In the \emph{Motion planning} example, using a grid size of 0.1 succeeds in synthesis while using 0.2 fails for abstraction-based methods because abstractions are more conservative for larger grid size. In contrast, our proposed method solves the problem in 151 seconds by using $\varepsilon=0.2$. This is because the minimum width of the partitions can be less than 0.2 by the bisection criterion in \cite[Algorithm 1]{Li2018rocs}. As opposed to the \emph{Motion planning} case where obstacles are distributed evenly across the state space, the constraints for parallel parking are highly nonlinear and only posed to a corner of the state space, and varying the discretization precision of the state space will save computational time in a great deal. Such a difference in those two case settings explains why the gain in time efficiency by using our method is more profound in the parallel parking cases.

As seen in (\ref{eq:complexity-abst}) and (\ref{eq:complexity-interval}), both methods are equivalently sensitive to the size of the discretized systems. The experimental results shows that the worst case as in (\ref{eq:complexity-interval}) is rather pessimistic in practice and our proposed method is more scalable to the discretization precision than abstraction-based methods. Analyzing the example of \emph{DC-DC converter}, we can observe in the right-hand side of Fig. \ref{fig:e.vs.t} that the run time of the proposed method changes slowly while the one for abstraction-based method explodes as precision $\varepsilon$ decreases. The left-hand side of Fig. \ref{fig:e.vs.t} compares the run time of the proposed method for two cases of different sizes, which indicates the dimensionality problem of the proposed method.

\begin{figure}[htbp]
  \centering
  \includegraphics[scale=0.45]{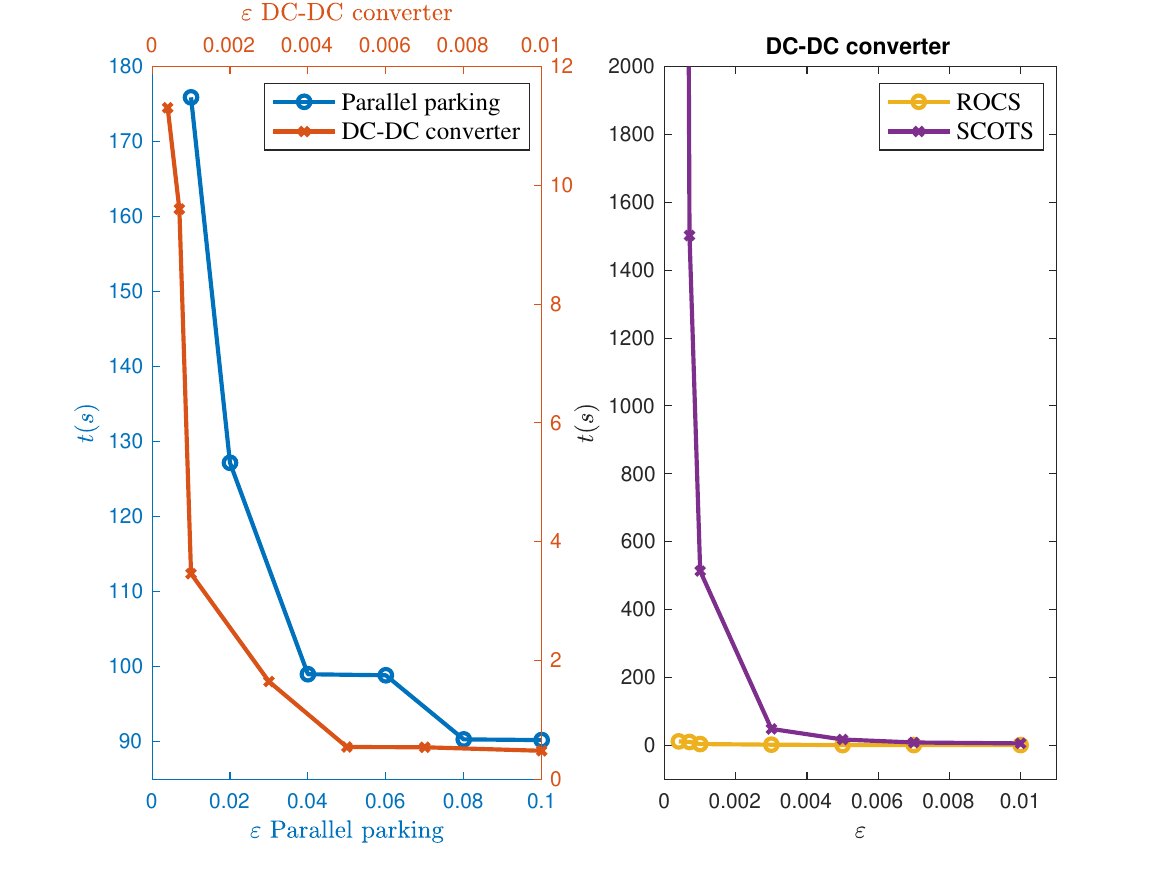}
  \caption{Changes of run time under different precisions.}
  \label{fig:e.vs.t}
\end{figure}
}

\section{Conclusions}
Under mild assumptions, we derived conditions so that reach-and-stay control synthesis is sound and robustly complete for discrete-time nonlinear systems in the sense that control strategies can be found if the specification can be satisfied for the perturbed system. A fixed-point algorithm based on interval computation was proposed as a practical control synthesis method. This is an improvement over abstraction-based methods, which are often not complete for systems without incremental stability. By adaptively partitioning the state space with respect to both dynamics and the given specification, the winning set for the given reach-and-stay problem can be inner-approximated with sufficiently high precision while reducing computational burdens. The efficiency was substantiated by performance tests on several benchmarking examples.

\bibliographystyle{IEEEtran}
\bibliography{../../reach}

%
\IEEEpeerreviewmaketitle









\end{document}